\documentclass[10pt]{article} 
\usepackage[preprint]{rlc}

\usepackage{amssymb}            
\usepackage{mathtools}          
\usepackage{mathrsfs}           
\mathtoolsset{showonlyrefs}     
\usepackage{graphicx}           
\usepackage{subcaption}         
\usepackage[space]{grffile}     
\usepackage{url}                
\usepackage[linesnumbered,ruled,vlined]{algorithm2e}

\title{Hybrid Reinforcement Learning Framework for Mixed-Variable Problems}

\author{Haoyan Zhai  \\
    haoyanzhai@gmail.com \\
    School of Mathematics\\
    Georgia Institute of Technology
    \And
    Qianli Hu \\
    qianlihuwork@gmail.com\\
    School of Mathematics\\
    Georgia Institute of Technology
    \And
    Jiangning Chen \\
    jiangning.chen@ukg.com\\
    UKG AI\\
    Ultimate Kronos Group}

\begin{document}

\maketitle

\begin{abstract}
Optimization problems characterized by both discrete and continuous variables are common across various disciplines, presenting unique challenges due to their complex solution landscapes and the difficulty of navigating mixed-variable spaces effectively. To Address these challenges, we introduce a hybrid Reinforcement Learning (RL) framework that synergizes RL for discrete variable selection with Bayesian Optimization for continuous variable adjustment. This framework stands out by its strategic integration of RL and continuous optimization techniques, enabling it to dynamically adapt to the problem's mixed-variable nature. By employing RL for exploring discrete decision spaces and Bayesian Optimization to refine continuous parameters, our approach not only demonstrates flexibility but also enhances optimization performance. Our experiments on synthetic functions and real-world machine learning hyperparameter tuning tasks reveal that our method consistently outperforms traditional RL, random search, and standalone Bayesian optimization in terms of effectiveness and efficiency.
\end{abstract}

\section{Introduction}

Optimization is a cornerstone of modern technological advancements, underpinning critical applications in machine learning \citep{sun2019survey}, operation research \citep{rardin1998optimization}, finance \citep{cornuejols2005optimization}, and engineering \citep{deb2012optimization}. As the problem complexity in these domains escalates - along with the exponential growth in the data sizes - the limitations of traditional optimization methods become increasingly apparent. This is particularly true for challenges that require tuning discrete and continuous variables simultaneously. While a plethora of optimization techniques have been developed to address either discrete or continuous variables effectively, seamlessly optimizing across a space imbued with both remains a formidable challenge.

Historically, researchers have treated continuous and discrete optimization tasks separately. On the one hand, Reinforcement Learning (RL) has carved its niche as a powerful tool for discrete optimization problems, lauded for its ability to learn and adapt strategies through direct interaction with the environment, see \citet{neunert2020continuous, mazyavkina2021reinforcement, poupart2006analytic}. Parallelly, continuous optimization techniques, notably Bayesian Optimization (BO), Gradient Descent, and Evolutionary Algorithms, have established their efficacy in navigating the subtleties of continuous parameter space \citep{andreasson2020introduction, zhan2022survey, omidvar2015designing, lauer2011continuous, kistler2003continuous}.

Despite these advancements, the quest for a unified approach that adeptly handles the intricacies of mixed-variable optimization tasks remains. These standalone applications of either discrete or continuous optimization methods tend to falter when faced with the multifaceted nature of real-world problems, underscoring a distinct gap in the current optimization toolbox \citep{liao2013ant, dimopoulos2007mixed, lin2018hybrid}.

\begin{figure}[ht]
  \centering
  \small
  \includegraphics[width=\linewidth]{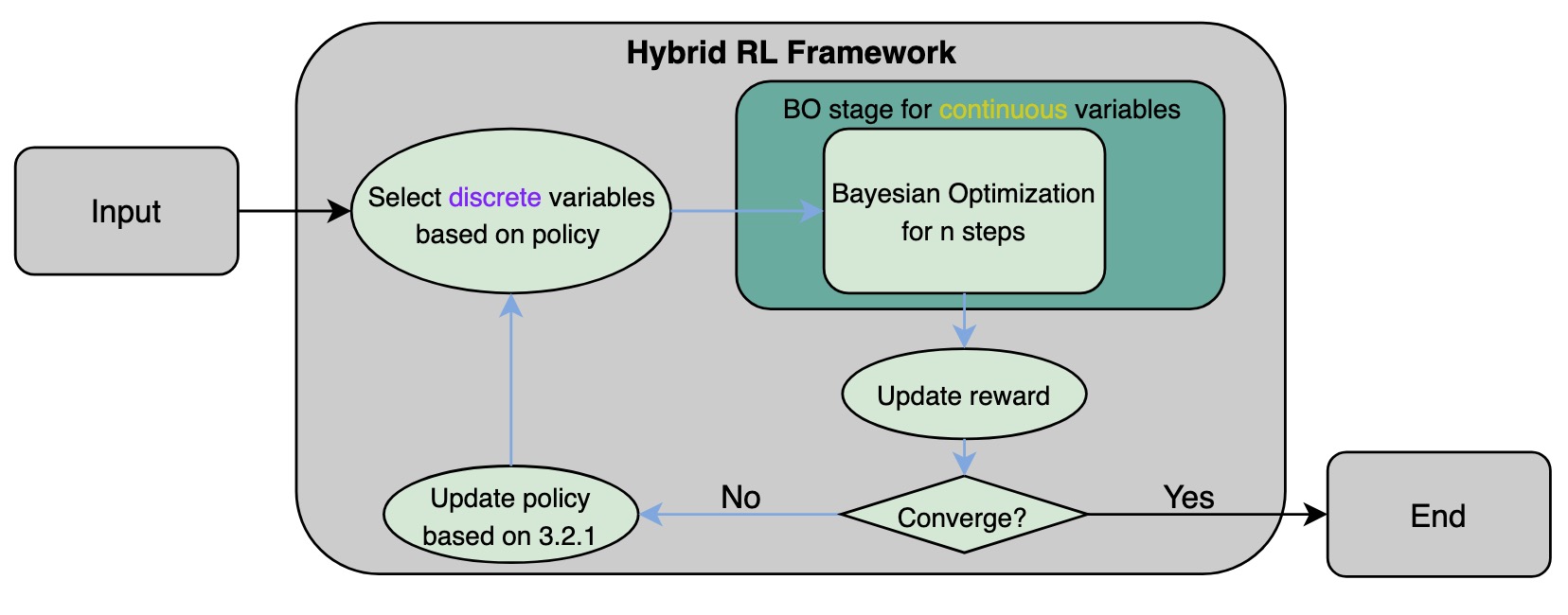}
    \caption{Overview of our hybrid Reinforcement Learning system. Each RL iteration is augmented by an n-step Bayesian Optimization for continuous variables until the stop criterion is met.}
  \label{fig:system_overview}
\end{figure}

To bridge this divide, we introduce a flexible hybrid Reinforcement Learning framework (depicted in Figure \ref{fig:system_overview}) that melds the strategic acumen of RL with the precision of continuous optimization techniques. This innovative approach not only navigates the complex terrain of mixed-variable problems with dexterity but also sheds some light on a potential paradigm in optimization, where the synergistic potential of integrating discrete and continuous methodologies is fully realized \citep{wang2001effective, piccoli1999necessary, cavazos2006hybrid, kelner2008hybrid}.

Our contributions through this work are manifold. Firstly, we introduce a novel Hybrid Framework that aims to integrate Reinforcement Learning (RL) with a diverse set of continuous optimization techniques, which provides a comprehensive strategy for effectively addressing mixed-variable optimization challenges. Secondly, our contributions are substantiated through empirical validation, involving rigorous experimentation on both synthetic functions and practical scenarios. The results consistently demonstrate the superior performance of our Hybrid Framework when compared to traditional RL methods, random search approaches, and standalone continuous optimization methods. Lastly, we highlight the adaptability and generalization of our framework across various continuous optimization strategies, showcasing its versatility and applicability across a wide range of optimization landscapes.

This paper is structured to guide the reader through the evolution of our framework and its validation. Section 2 provides a review of related work, setting the stage for our methodology. Section 3 delves into the flexible hybrid Reinforcement Learning framework, detailing its components and operational dynamics. Section 4 outlines our experimental design and the ensuing results, highlighting the effectiveness and versatility of the approach. Section 5 concludes the paper, summarizing our contributions and envisioning the path forward.

\section{Related Work}

The exploration of optimization strategies within artificial intelligence has led to significant advancements across various domains. Reinforcement Learning (RL) and continuous optimization techniques stand out as pivotal components in this evolutionary trajectory. This section reviews the pertinent literature, focusing on the development and application of RL in optimization, the role of continuous optimization methods, and the inception of hybrid frameworks that seek to leverage the strengths of both approaches.

Reinforcement Learning has emerged from its roots in dynamic programming to become a key methodology for solving decision-making problems characterized by uncertainty and complexity. The foundational work by \citet{sutton2018reinforcement} provides a comprehensive overview of RL principles, outlining how agents learn to make decisions through interactions with their environment to maximize some notion of cumulative reward.

In the context of optimization, RL has been particularly impactful in discrete spaces, where the goal is to find optimal sequences or configurations from a finite set of possibilities. Applications range from combinatorial optimization \citep{mazyavkina2021reinforcement} to more nuanced domains like network configuration \citep{gao2019dynamic, mammeri2019reinforcement} and inventory management \citep{giannoccaro2002inventory, sultana2020reinforcement}. These studies underscore RL's versatility and capacity to adapt strategies based on feedback, a trait that is invaluable for navigating the complex landscapes of optimization problems.

Parallel to the advancements in RL, continuous optimization techniques have been refined and adapted to meet the challenges of tuning continuous parameters. Techniques such as Bayesian Optimization (BO), Gradient Descent, and Evolutionary Algorithms have been at the forefront of this endeavor. Bayesian Optimization, in particular, has gained prominence for its efficiency in hyperparameter tuning of machine learning models \citep{wu2019hyperparameter}, leveraging probabilistic models to guide the search process. Gradient Descent remains a staple in training neural networks by iteratively minimizing loss functions \citep{bengio2000gradient}, while Evolutionary Algorithms offer a population-based search strategy, drawing inspiration from biological evolution to optimize complex functions \citep{young2015optimizing, tani2021evolutionary}.

Recognizing the complementary strengths of RL and continuous optimization techniques, recent research has begun to explore hybrid frameworks that integrate these methodologies. Such frameworks aim to capitalize on RL's ability to make strategic decisions in discrete spaces, while employing continuous optimization methods for fine-tuning continuous variables. This synergy has the potential to offer more robust and efficient solutions to mixed-variable optimization problems, a frontier that is increasingly relevant in areas such as automated machine learning \citep{liu2020admm} and dynamic resource allocation \citep{bagci2018dynamic}.

Initial attempts at such integration have shown promising results, suggesting that the combined approach could outperform the application of either methodology in isolation \citep{lin2018hybrid, wang2017hybrid, wang2021particle}. These hybrid models not only offer a nuanced understanding of the optimization landscape but also open new avenues for research in both theoretical development and practical applications.

Despite these advancements, there remains a notable gap in the literature concerning efficient, scalable, and adaptable frameworks that can effectively navigate the complexities of mixed-variable optimization landscapes. The development of a flexible hybrid Reinforcement Learning framework, as proposed in this paper, aims to address this gap by offering a novel integration of RL and continuous optimization techniques, tailored to the demands of contemporary optimization challenges.

\section{Methodology}

\subsection{Problem Definition}

This study focuses on optimization problems that feature a combination of discrete and continuous decision variables, with the aim of maximizing an objective function. Let $\mathcal{A}$ represent the space of discrete variables and $\mathcal{X}$ the space of continuous variables. The objective is to find the optimal pair $(a^*, x^*)$ that maximizes the objective function $f: \mathcal{A} \times \mathcal{X} \rightarrow \mathbb{R}$, which can be formally expressed as:

\[
\max_{a \in \mathcal{A}, x \in \mathcal{X}} f(a, x)
\]

\subsection{Hybrid Reinforcement Learning Framework}

The flexible hybrid Reinforcement Learning framework addresses mixed-variable optimization problems by bifurcating the optimization process: it employs Reinforcement Learning (RL) for discrete variable selection and a specific continuous optimization technique for continuous variable adjustment. Herein, we employ the Gradient Bandit model with Softmax action selection for RL, and Bayesian Optimization for updating the reward function.

\subsubsection{Reinforcement Learning for Discrete Variable Optimization}

For the discrete variable optimization component, we utilize the Gradient Bandit agent, which is designed to select actions that optimize cumulative rewards over time. The action selection is governed by the Softmax policy, defined as follows:

\begin{equation}
    Pr\{A_t=a\}\frac{\exp(H_t(a))}{\sum_{b=1}^n\exp(H_t(b))}=\pi_t(a),
\label{equ:pi}
\end{equation}

where $A_t$ represents the action taken at time $t$, $\pi_t(a)$ the probability of taking action $a$ at time $t$, and $H_t(a)$ the learned numerical preference for each action $a$. The preferences are updated upon selecting action $A_t$ and observing reward $R_t$, using the formulae:

\begin{align}
    &H_{t+1}(A_t)=H_t(A_t)+\alpha(R_t-\Bar{R_t})(1-\pi_t(A_t)),&and\\
    &H_{t+1}(a)=H_t(a)-\alpha(R_t-\Bar{R_t})\pi_t(a),&\forall a\neq A_t
\label{equ:H}
\end{align}

where $\alpha > 0$ is a step-size parameter, and $\Bar{R}_t\in\mathbb{R}$ denotes the average of all rewards up to and including time $t$ \citep{sutton1999reinforcement}. This mechanism ensures an exploration-exploitation balance in action selection, enabling the agent to prioritize actions based on their reward potential. Through interactions with the objective function, the RL agent dynamically refines its policy based on reward feedback, enhancing its decision-making accuracy over time.

\subsubsection{Continuous Optimization Techniques}

A crucial aspect of the Gradient Bandit's effectiveness is defining the reward function $R_t$. For this, we employ Bayesian Optimization. After setting the discrete variables $A_t$, we conduct several iterations of Bayesian Optimization on the problem: $\max_{x\in\mathcal{X}}f(A_t,x)$, retaining the solver status for $A_t$. If $A_t$ has been previously encountered, we reload its continuous problem status and continue optimization; otherwise, we initiate a new problem. The outcome of Bayesian Optimization, denoted as $\{f(A_t,X_1), f(A_t,X_2), \ldots, f(A_t,X_s)\}$, informs the definition of the reward:

\begin{equation}
    R_t(A_t)=\max_sf(A_t,X_s).
    \label{equ:reward}
\end{equation}

The whole algorithm is given as Alg.\ref{alg:hybrid}. While this methodology employs the Gradient Bandit and Bayesian Optimization for demonstration, it is important to note that the framework is flexible and can incorporate other RL methods and continuous optimization algorithms to suit different problem contexts and requirements.

\begin{algorithm}[ht]
\SetAlgoLined
\caption{Hybrid RL Optimizer}
\label{alg:hybrid}
\small
\DontPrintSemicolon
\KwIn{\\
    $f(a,x)$,\\
    Bayesian Optimization Iteration Steps $n$,\\
    $H_0(a)=0$,\\
    $Cache[a]=\{\}$, $\forall a\in\mathcal{A}$,
}
\For{$t=0,1,2,3\cdots$}{
    Calculate $\pi_t(a)$ $\forall a\in\mathcal{A}$.\\
    Choose an action $A_t\in\mathcal{A}$ according to the probability distribution $\pi_t$.\\
    \eIf{$Cache[A_t]$ is empty}{
        Initialize a Bayesian Optimizer $opt$.
    }{
        $opt=Cache[A_t]$.
    }
    Iterate $opt$ for $n$ steps, and get the reward $R_t=\max_{s\leq t}f(A_t,X_{sn})$, where $X_{sn}$ is the iteration results in $opt$.\\
    Save the current Bayesian Optimization status $Cache[A_t]=opt$.\\
    Update $H_t$ according to \ref{equ:H}.\\
    \If{$(X_t,R_t(X_t))$ appears $m$ times in the last $T$ iterations}{
        Break
    }
}
\end{algorithm}

\section{Experiments}

\subsection{Experimental Setup}


The experiments were conducted on a computing environment equipped with AMD Ryzen Threadripper 3970X 32-Core Processor with 256GB RAM. The framework was implemented in Python, with reinforcement learning components based on \citet{galbraith_2023}, and continuous optimization techniques, including Bayesian Optimization, Gradient Descent, and Evolutionary Algorithms, were utilized as per the framework's design.

\subsection{Datasets and Scenarios}

To assess the framework's performance, we use a combination of synthetic functions and real-world scenarios. In the experiments, we use Random Search, Gradient Bandit, Bayesian Optimization, and our Hybrid method to perform the optimization problems. Since Gradient Bandit only takes discrete variables, we discretize continuous variables in their corresponding domains. For Bayesian Optimization, we choose variable values in the range of the minimum and the maximum of the variable, and round it to the nearest discrete value in the domains.

\subsubsection{Synthetic Functions} Three Benchmark functions are employed to evaluate the framework's capability to find global optima in complex landscapes. The objective function is simply the function values, and we report the absolute difference between the calculated function values and the known global maximum in Section \ref{sec:results}

\begin{enumerate}
\item \textbf{Shekel Function:} The first one is the Shekel function, which is a benchmark test function in global optimization problems. The function is given as below:
\begin{equation}
    f_1(x_1,x_2,x_3,x_4)=\sum_{i=1}^4\left(c_i+\sum_{j=1}^4(x_j-a_{ji})^2\right)^{-1}, 
\label{shekel}
\end{equation}
where
\[
C=(c_1,\cdots,c_{10})=\frac{1}{10}(1,2,2,4,4,6,3,7,5,5),
\]
and
\[
A=(a_{ji})=\left(
\begin{array}{cccccccccc}
    4 & 1 & 8 & 6 & 3 & 2 & 5 & 8 & 6 & 7\\
    4 & 1 & 8 & 6 & 7 & 9 & 3 & 1 & 2 & 3.6\\
    4 & 1 & 8 & 6 & 3 & 2 & 5 & 8 & 6 & 7\\
    4 & 1 & 8 & 6 & 7 & 9 & 3 & 1 & 2 & 3.6\\
\end{array}\right).
\]
To make it a function containing both discrete and continuous variables, we set the first two variables to be integers. This function has 10 local maximums and the global maximum is around $10.536283726219603$ at $(4,4,4,4)$, and the searching region is $\{0,1,2,3,4,5,6,7,8,9,10\}\times\{0,1,2,3,4,5,6,7,8,9,10\}\times[0,10]\times[0,10]$.

\item \textbf{Composition Function:} The second function is a composition of Rastrigin (global minimum $f(0,0)=0$), Ackley (global minimum $f(0,0)=0$), and Sphere (global minimum $f(0,0)=0$) functions:
\begin{align}
&\text{Rastrigin}(\mathbf{x})=10 n + \sum_{i=1}^n \left[x_i^2 - 10\cos(2 \pi x_i)\right],\ \ \ \text{Sphere}(x,y)=x^2+y^2\\
&\text{Ackley}(x,y)=-20\exp\left[-0.2\sqrt{0.5(x^2+y^2)}\,\right] -\exp\left[0.5\left(\cos 2\pi x + \cos 2\pi y \right)\right] + e + 20\\
\end{align}
We define our function as:
\begin{equation}
f_2(u,x,y)=\left\lbrace\begin{array}{ll}
-\text{Rastrigin}(x,y) & u=0 \\
-\text{Ackley}(x,y)+10 & u=1 \\
-\text{Sphere}(x,y)+20 & u=2 \\
\end{array}\right.
\label{mix}
\end{equation}

Therefore, the global maximum is $f(2,0,0)=20$. We also define the variable $x$ to be discrete and the search area to be $\{0,1,2\}\times\{-1,0,1,2,3\}\times[-5,5]$.

\item \textbf{Sine Permutation Function:} At last, we define a two-variable function as below
\[
g(x,y)=\frac{x\sin\left((-x+7)^2\pi/(2(y-4)^2+1)\right)}{((x-5)^2+1)}.
\]
Then, we define $u,v,w\in\{1,4,7,10,13\}$, $x\in[0.5,8]$, and $y\in[0.1,5]$. Define the permutations $Perm_u(x)=(7,1,13,10,4)$, in which the value of this function is the integer next to the given integer $x$ in the array, i.e. $Perm_u(1)=13$. Following the same logic, also define $Perm_v(x)=(13,1,4,7,10)$, and $Perm_w(x)=(7,4,10,1,13)$. Now we can define our objective function as
\begin{equation}
    f_3(u,v,w,x,y)=g(Perm_u(u)+Perm_v(v)+Perm_w(w)+x,y).
\end{equation}
\end{enumerate}

\subsubsection{Real-World Scenario} Referred to the Walmart dataset in Kaggle \citep{yasserh_walmart_dataset}, we predict the upcoming weekly sales with the tree model XGBoost. Seeking suitable hyperparameters to minimize the mean absolute error (MAE) of the regression problem, we tune $max\_bin$, $n\_estimators$, and $max\_depth$, which are discrete parameters, incorporating continuous parameters $colsample\_bytree$, $learning\_rate$, $reg\_alpha$, and $reg\_lambda$. To prevent over-fitting and under-fitting, we split the whole dataset as a training set and a validation set. Fixing a set of hyperparameters, we train the model on the training dataset and then calculate the MAE of the validation set. As our framework is maximizing the objective and to avoid number overflow, we set the objective function of this problem as $-\text{MAE}/1e5$.

\subsection{Results}
\label{sec:results}

\begin{figure}[ht]
  \centering
  \includegraphics[width=\linewidth]{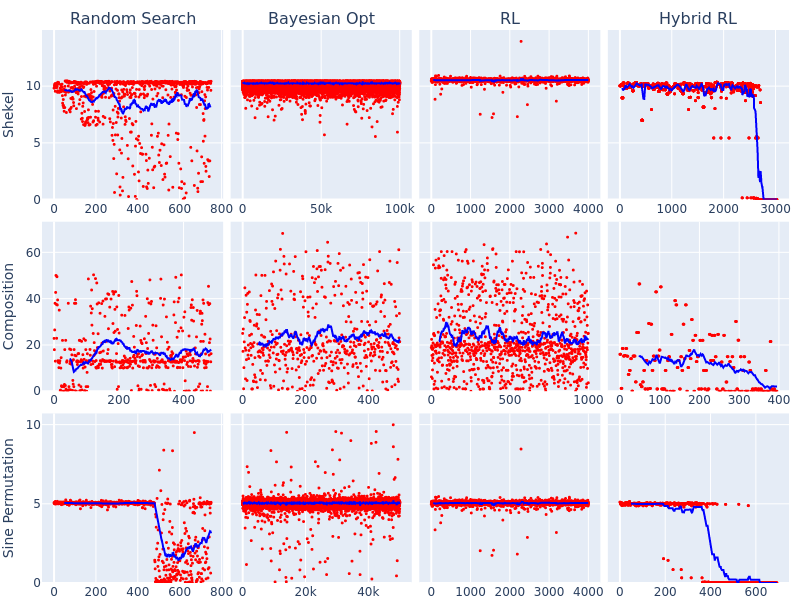}
    \caption{The results for synthetic functions. Each graph is a trajectory of the optimization. The x-axis is the number of iterations, while the y-axis is the absolute gap between the searched maximum and the known global maximum. The red dots are the result of one step of the experiment, and the blue curve is the rolling average of the adjacent 50 points. Specifically in hybrid RL, we need to define the Bayesian Opt steps to perform in each reinforcement learning step $n$: for Shekel function $n=3$, for Composition Function $n=3$, for Sine Permutation Function $n=2$.}
  \label{fig:synthetic_one_exp}
\end{figure}

\begin{figure}[ht]
  \centering
  \includegraphics[width=\linewidth]{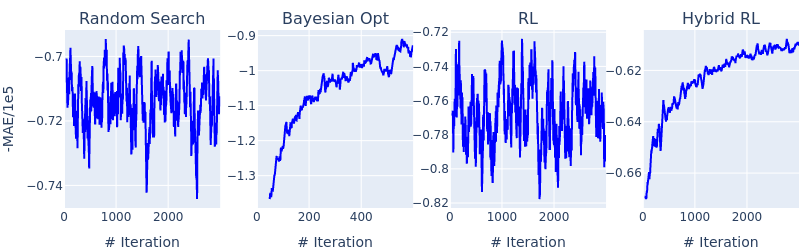}
    \caption{The results for machine learning hyperparameter tuning problem. Each curve is the rolling average of the adjacent 50 points. The x-axis is the number of iterations, while the y-axis is the reward. The reward is defined as the objective function value at each step for Random Search, Bayesian Opt, and RL, while it is defined as Equation \ref{equ:reward}. In Hybrid RL, the Bayesian Opt step $n=2$.}
  \label{fig:ml_exp}
\end{figure}

\begin{figure}[ht]
  \centering
  \includegraphics[width=\linewidth]{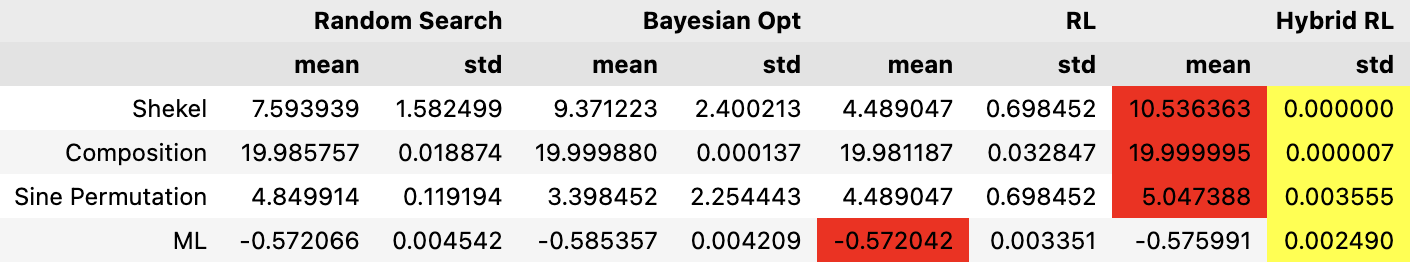}
    \caption{The statistic result for all the problems repeated 5 times for synthetic functions and 10 times for the machine learning problem with different random seeds. We can see that the hybrid RL gives the lowest std in all problems and achieves the optimal values in most cases.}
  \label{fig:repeated_exp_stats}
\end{figure}

We have conducted four sets of optimization experiments detailed in section 4.2.1, covering the Shekel function, a composition function, the Sine Permutation function, and a machine-learning-based hyper-parameter tuning task). Each experiment is evaluated by the same four optimization methods: random search, Bayesian Optimization, vanilla Reinforcement Learning, and our hybrid method. The results demonstrate that our flexible hybrid Reinforcement Learning framework significantly outperforms the comparative methods. More specifically:

\textbf{Convergence Assurance:} {As shown in Figure \ref{fig:synthetic_one_exp}, our hybrid RL method demonstrates a high level of convergence certainty. The experiments support the assertion that the hybrid method has amalgamated the benefits of reinforcement learning and conventional optimization methods tailored for continuous variables such as Bayesian Optimization.}


\textbf{Convergence Efficiency:} {Our hybrid method, enhanced with periodic Reinforcement Learning checks, outpaces Bayesian Optimization in computational speed despite requiring more iterations, as Bayesian Optimization's per-iteration time increases with expanded search spaces, evident in Figure \ref{fig:ml_exp}. This results in our hybrid approach achieving faster convergence within practical time frames, even beyond 600 rounds where Bayesian Optimization becomes untenably slow. Compared to vanilla Reinforcement Learning, which necessitates segmenting continuous spaces into discrete intervals—thus expanding the action space and computational load—our hybrid method maintains the action space's continuity. This unique balance allows it to converge in fewer iterations than pure RL, despite the individual iterations taking slightly longer, showcasing an optimal trade-off between iteration count and time efficiency for complex optimization tasks.}

\textbf{Low Standard Deviation Between Rounds of Experiments:} {Compared to conventional methods, our hybrid approach demonstrates remarkable consistency, reflected in its consistently low standard deviations across the four experiments. Figure \ref{fig:repeated_exp_stats} illustrates that although our hybrid model lags slightly behind in the ML hyper-parameter tuning task, the statistical analysis of standard deviations suggests a notably stable performance. Furthermore, in the remaining three tests, the hybrid model consistently identifies the global optimum.}



\section{Conclusion and Future Research}

This paper presented a flexible hybrid Reinforcement Learning Framework that innovatively combines reinforcement learning with continuous optimization techniques such as Bayesian Optimization to tackle mixed-variable optimization problems. Through our experimentation, we have demonstrated its efficacy and efficiency over traditional approaches. The framework's flexibility and adaptability suggest its potential applicability across a broad spectrum of optimization challenges. We believe that this work not only advances the fusion of continuous optimization techniques and reinforcement learning but also opens new doors for solving complex problems in various domains. Future research will aim to expand the framework's capabilities and its application in more diverse and large-scale scenarios, to theoretically substantiate the efficacy of our algorithm, and to explore more sophisticated reinforcement learning and continuous optimization algorithms to further improve the framework's effectiveness.

\bibliography{main}
\bibliographystyle{rlc}






\end{document}